\documentclass{IEEEtran}
\usepackage{cite}
\usepackage{amsmath,amssymb,amsfonts}
\usepackage{algorithmic}
\usepackage{graphicx}
\usepackage{textcomp}
\usepackage{hyperref, comment}
\usepackage{textcomp}
\usepackage{xcolor}
\usepackage{dsfont}

\newtheorem{lemma}{Lemma}
\newtheorem{theorem}{Theorem}

\newtheorem{remark}{Remark}
\newtheorem{definition}{Definition}

\newtheorem{problem}{Problem}
\newtheorem{condition}{Condition}

\usepackage{dblfloatfix}
\def\BibTeX{{\rm B\kern-.05em{\sc i\kern-.025em b}\kern-.08em
    T\kern-.1667em\lower.7ex\hbox{E}\kern-.125emX}}
\begin{document}

\title{Robust Output Regulation with \\  Ellipsoidal State and Input Constraints}


\author{K\"{u}r\c{s}ad Metehan G\"{u}l and Selahattin Burak Sars\i lmaz$^{}$
\thanks{$^{}$ K\"{u}r\c{s}ad Metehan G\"{u}l and Selahattin Burak Sars\i lmaz are with the Department of Electrical and Computer Engineering, Utah State University, Logan, UT 84322 USA (emails: {\tt kursad.gul@usu.edu and burak.sarsilmaz@usu.edu).}}
}

\maketitle

\begin{abstract}
This letter studies robust output regulation of both continuous-time and discrete-time uncertain linear time-invariant systems subject to ellipsoidal state and input constraints. We show that the existence of an open set of uncertainties containing the zero vector, in which the closed-loop system is internally stable and an augmented ellipsoidal set for the closed-loop system is robustly positively invariant, together with an inequality for input constraint satisfaction, are sufficient for the solvability of the problem with an internal model-based control law. Here, the augmented ellipsoidal set is obtained by augmenting the system's state constraint set with a virtual state constraint set of the internal model. 
For both continuous-time and discrete-time cases, we derive linear matrix inequalities that ensure the existence of such uncertainty sets and directly lead to convex semidefinite feasibility programs for controller parameter synthesis, which can be reliably solved using off-the-shelf convex solvers. A numerical example demonstrates a selected result of this letter and the corresponding control parameter synthesis.

\end{abstract}



\section{Introduction}
\subsection{Motivation and Literature Review}

In practice, automatic control of dynamical systems involves satisfying state and input constraints imposed by workspace and actuator limitations \cite{application_spacecraft,application_quadrotor_mobilerobot}. To address these constraints while guaranteeing closed-loop stability, researchers often resort to the concept of positive invariance \cite{Blanchinisurvey,Blanchinisbook,RAKOVIC2007831}. Fundamental studies \cite{1104270,VASSILAKI01061988} leveraged this concept to 
address control of discrete-time linear systems with polyhedral state or input constraints.  They were later extended to consider potentially unmatched additive disturbances  \cite{Blanchini_constrained_1} and uncertain system and input matrices \cite{MILANI19951489}. More recently, \cite{9763329} revealed necessary and sufficient conditions for robust positive invariance of a polyhedral set of the state based on noisy input-state data. 

Ellipsoidal sets are also of particular interest as candidate positively invariant sets since they are associated with Lyapunov equations and linear matrix inequalities (LMIs), which are convex semidefinite feasibility programs \cite{BoydLmibook,Kurzhanski_ellipsoidal}. In \cite{pordal2025ellipsoidalsettheoreticdesignrobust}, an LMI formulation for robust controlled positively invariant set synthesis for continuous-time linear systems has been given. Inscribed positively invariant ellipsoidal sets have been frequently used to address control problems with state and input constraints, since they can be computed through convex programs \cite {KOTHARE19961361,BROOMS2001157,pordal2025ellipsoidalsettheoreticdesignrobust,Llorente_BLF}.

Robust output regulation \cite{1101137, FRANCIS1976457, IsidoriORP, huang2004nonlinear, 10502160} offers a unified framework for formulating the problem of asymptotic tracking and rejection of a broad class of signals with closed-loop stability, while accounting for the effect of parametric system uncertainties on steady-state tracking performance. The output regulation problem (ORP), which aims to design a control law such that the closed-loop system is stable and the system’s output asymptotically tracks a class of reference signals while rejecting a class of disturbances, has been mainly treated via feedforward and internal model approaches. The feedforward approach \cite{IsidoriORP,huang2004nonlinear} employs a feedforward gain designed based on a simultaneous solution of a pair of linear matrix equations associated with the nominal system matrices, known as the regulator equations.
Consequently, the resulting control law does not exhibit inherent robustness to parametric uncertainties;  that is, even very small parametric uncertainties can lead to nonzero steady-state error although closed-loop stability is preserved. On the other hand, the internal model approach \cite{1101137, FRANCIS1976457, IsidoriORP, huang2004nonlinear, 10502160} augments the system dynamics with a dynamic control law that incorporates a model of references and disturbances (as in the celebrated proportional-integral (PI) controller incorporating the model of constant signals for perfect tracking and rejection of constant signals at steady state). This ensures robustness against parametric uncertainties: for parametric uncertainties that do not destabilize the closed-loop system, the steady-state tracking error remains zero.

Using the feedforward approach, the ORP of discrete-time linear time-invariant (LTI) systems with polytopic state and input constraints is addressed in \cite{DD_Constrainted_ORP_Chinese_Paper, WANG2023111050}. In \cite{10670461}, general compact and convex sets containing the origin in their interiors are allowed as state and input constraint sets by combining the model predictive control framework with the feedforward approach. Due to the feedforward approach, these studies may no longer achieve zero steady-state tracking errors even in the presence of small parametric uncertainties that do not destabilize the closed-loop system. 

A PI-like output feedback controller is adopted to address tracking of constant signals under polytopic state and input constraints for continuous-time LTI systems in \cite{10056979} and discrete-time LTI systems in \cite{dos2023discrete}. The method is then extended to track sinusoidal signals for single-output continuous-time LTI systems under polytopic state and input constraints \cite{10644900}. Although the underlying approach is the internal model approach, the studies  \cite{10056979, dos2023discrete, 10644900} do not consider parametric system uncertainties, disturbances, and direct feedthrough terms. 
The paper \cite{slidingmode_CRORP_SISO} addresses constant signal tracking for single-input single-output continuous-time LTI systems with bounded parametric uncertainty only in the system matrix and matched disturbances using 
an integral sliding-mode-based control law.
Despite these efforts, the literature lacks a unified formulation and a solution of the robust output regulation problem with the internal model approach for multi-input multi-output uncertain LTI systems subject to ellipsoidal state and input constraints.

\subsection{Contribution}
We provide a unified formulation and a systematic treatment of the constrained robust output regulation problem (CRORP) with the internal model approach for both continuous-time and discrete-time general uncertain LTI systems under ellipsoidal state and input constraints. This formulation allows  tracking and rejection of a broad class of references and disturbances, including constant signals with any amplitudes, sinusoidal signals with different amplitudes, frequencies, and phases, and any linear combinations of these signals.
While model predictive control, commonly employed to handle constrained control problems, typically requires repeated online optimization via a receding-horizon implementation, the results of this letter provide a one-shot, offline controller synthesis that eliminates the need for repeated online optimization.






We first derive a sufficient condition for the existence of an open set of uncertainties containing $0$, in which an augmented ellipsoidal set is robustly positively invariant for the closed-loop system (see Lemma \ref{lmm:ct_strict_invariance} for the continuous-time case and Lemma \ref{lmm:dt_strict_invariance} for the discrete-time case). Here, the augmented ellipsoidal set is obtained by augmenting the system's state constraint set with a virtual state constraint set of the internal model.  Building on this result, we then provide an LMI condition that ensures not only the existence of such an uncertainty set but also the internal stability of the closed-loop system within this uncertainty set (see Lemma \ref{lmm:ct_time_RPI} for the continuous-time case and Lemma \ref{lmm:dt_time_RPI} for the discrete-time case). 
 Our main result establishes that this LMI, together with another LMI for input constraint satisfaction, is sufficient for the solvability of the CRORP under standard assumptions (see Theorem \ref{thm:main_thm_ct_time} for the continuous-time case and Theorem  \ref{thm:main_thm_dt_time} for the discrete-time case).  The derived sufficient LMI conditions directly lead to convex semidefinite feasibility programs for controller parameter synthesis for both continuous-time and discrete-time cases, which can be reliably solved with off-the-shelf convex solvers. Finally, we provide a numerical example to demonstrate a selected result of this letter and the corresponding control parameter synthesis.

\subsection{Notation and Definitions}
We denote the set of nonnegative real numbers with $\mathbb{R}_{+}$, the set of nonnegative integers with $\mathbb{Z}_{+}$, and the set of $n \times n$ symmetric positive definite matrices with $\mathbb{S}^n_{++}$. We write $I_n$ for the $n \times n$ identity matrix,
and $\mathrm{diag} (X_1, \ldots, X_n)$ for  
the block diagonal matrix with matrices $X_1, \ldots, X_n$ on its diagonal blocks. For symmetric matrices $X$ and $Y$,  $X  \succ Y $ $(\succeq, \prec, \preceq)$ means $X-Y$ is positive definite (positive semidefinite, negative definite, negative semidefinite).  We write  $\lambda_{\mathrm{min}}(\cdot)$ and  $\lambda_{\mathrm{max}}(\cdot)$ for the minimum and maximum eigenvalue of a symmetric matrix. A square matrix $X \in \mathbb{R}^{n \times n}$ is said to be stable if it is Hurwitz\footnote{That is, all its eigenvalues lie in the open left half complex plane.} for continuous-time systems and Schur\footnote{That is, all its eigenvalues have modulus less than $1$.} for discrete-time systems. 
The vector formed by stacking columns of a matrix $X \in \mathbb{R}^{n \times m}$ is denoted by $\mathrm{vec}(X)$.  We denote the Euclidean norm for vectors and the induced $2$-norm for matrices by $\|\cdot\|$. The interior and boundary of a set $\Omega \subseteq \mathbb{R}^n$ are denoted by $\mathrm{int}\Omega$ and $\partial\Omega$, respectively. Finally, $\rho$ is an operator indicating the time derivative for continuous-time systems ($\rho x(t) = \dot{x}(t)$ for $t \in \mathbb{R}_+$) and the forward unit time shift for discrete-time systems ($\rho x(t) = x(t+1)$ for $t \in \mathbb{Z}_+$).

\begin{definition}[Definition 2.2 in \cite{Blanchinisurvey}]\label{def:RPIset}
    A set $\mathcal{X} \subseteq \mathbb{R}^n$ is called robustly positively invariant for the system
    \begin{equation}\nonumber
        \rho x = Ax + Bv        
    \end{equation}
    and the set $\mathcal{V} \subseteq \mathbb{R}^q$ provided that $x(t) \in \mathcal{X}$ for all $t>0$ whenever $x(0) \in \mathcal{X}$ and $v(t) \in \mathcal{V}$ for all $t \geq 0$, where $v$ is assumed to be continuous on $\mathbb{R}_{+}$ for continuous-time case. 
\end{definition}

\section{Problem Formulation}
Consider an uncertain system of the form
\begin{align} \label{eq:uncertain_system}\nonumber
\rho x&= {A}_w x + {B}_w u + {E}_w v \\
e &= {C}_w x + {D}_w u + {F}_w v, 
\end{align}
where $x \in \mathbb{R}^n$ is the state, $u \in \mathbb{R}^m$ is the control input, and $e \in \mathbb{R}^p$ is the tracking error. 
The matrices ${A}_w = A + \Delta A$, ${B}_w = B + \Delta B$, ${C}_w = C + \Delta C$, ${D}_w = D + \Delta D$, ${E}_w = E + \Delta E$, and ${F}_w = F + \Delta F$ are uncertain matrices with $A$, $B$, $C$, $D$, $E$, and $F$ being their nominal parts.  Also, 
$v \in \mathbb{R}^q$ is the exogenous signal generated by the exosystem
\begin{equation} \label{eq:exosystem}
\rho v = A_0 v
\end{equation}
and it is assumed to lie entirely in a known ellipsoidal set 
\begin{align}\nonumber
    \mathcal{V} &= \{v \in \mathbb{R}^{q} \ | \ v^\mathrm{T}R v \leq 1 \}, \ R \in \mathbb{S}^{q}_{++}, 
\end{align}
i.e., $v(t) \in \mathcal{V}$ for all $t \geq 0$. The exosystem 
yields the reference $-{F}_wv$ to be tracked and the disturbance  ${E}_wv$ to be rejected by the uncertain system \eqref{eq:uncertain_system}.
 For convenience, the system uncertainties  are represented by 
    $w = [
        \mathrm{vec}(\Delta A)^\mathrm{T}, \mathrm{vec}(\Delta B)^\mathrm{T}, \mathrm{vec}(\Delta C)^\mathrm{T}, \mathrm{vec}(\Delta D)^\mathrm{T}, \mathrm{vec}(\Delta E)^\mathrm{T}, \\ \mathrm{vec}(\Delta F)^\mathrm{T}
    ]^\mathrm{T} \in \mathbb{R}^{\bar n}$, where 
$\bar n = n(n+m+p+q)+p(m+q)$.

The ellipsoidal state and input constraint sets considered in this study are as follows:
\begin{align}\nonumber
    \mathcal{X} &= \{x \in \mathbb{R}^n \ | \ x^\mathrm{T}\Gamma x \leq 1 \}, \ \Gamma \in \mathbb{S}^n_{++} \\
    \mathcal{U} &= \{u \in \mathbb{R}^m \ | \  u^\mathrm{T}\Pi u \leq 1 \}, \ \Pi \in \mathbb{S}^m_{++}.  \nonumber
\end{align}
We consider the following dynamic state-feedback control law 
\begin{align}\label{eq:controllaw} \nonumber
    \rho z &= G_{1}z + G_{2}e  \\
    u &= K_{1}x + K_{2}z + K_3v, 
\end{align}
where $z \in \mathbb{R}^{n_z}$ is the controller state. The matrices $G_{1} \in \mathbb{R}^{n_z \times n_z}$,  $G_{2} \in \mathbb{R}^{n_z \times p}$, $K_{1} \in \mathbb{R}^{m \times n}$, $K_{2} \in \mathbb{R}^{m \times n_z}$, and $K_{3} \in \mathbb{R}^{m \times q}$ are  control gains to be designed\footnote{The control law \eqref{eq:controllaw} includes a feedforward term $K_3v$, which requires knowledge of the exosystem state $v$. However, this term is optional, and the proposed results accommodate the case $K_3=0$. In this case, \eqref{eq:controllaw} reduces to the internal model-based dynamic state-feedback control law commonly used in the robust output regulation theory (e.g., see \cite{huang2004nonlinear,10502160}), eliminating the need for direct access to the exosystem state $v$. The inclusion of $K_3v$ provides additional design flexibility, particularly for improving transient performance, when the exosystem state $v$ is available (see Example 8.5 in \cite{ogata2010modern}).}.

We now formally state the CRORP. 

\begin{problem}[CRORP]\label{prb:CRORP}
    Given the system \eqref{eq:uncertain_system} and the exosystem \eqref{eq:exosystem}, find a dynamic state-feedback control law of the form \eqref{eq:controllaw} such that there exists an (open) neighborhood\footnote{In this letter, we adopt the definition of a neighborhood of a point given in Definition 2.18(a) of \cite{babyRudin}.} $\mathcal{W}$ of $w = 0$ for which, for any $x(0) \in \mathcal{X}$, $z(0) = 0$, $v(t) \in \mathcal{V}$ for $t\geq0$, and $w \in \mathcal{W}$, we have
    \begin{enumerate}
        \item [(i)] $x(t) \in \mathcal{X} \ \text{for all} \ t > 0, \  u(t) \in \mathcal{U} \ \text{for all} \ t \geq 0$, and
        \item [(ii)] $\lim_{t\to\infty} e(t) = 0$.
    \end{enumerate}
\end{problem}

Letting $\zeta = [x^\mathrm{T},z^\mathrm{T}]^\mathrm{T}$ and using \eqref{eq:uncertain_system} and \eqref{eq:controllaw}, one can obtain the closed-loop system 
\begin{align} \label{eq:augmented_clloop_system} 
\rho \zeta&= {A}_{\mathrm{c}w} \zeta + {B}_{\mathrm{c}w} v,
\end{align}
and the tracking error can be written as
\begin{align}\label{eq:error_initial}
e &= {C}_{\mathrm{c}w} \zeta + {D}_{\mathrm{c}w} v,
\end{align}
where ${A}_{\mathrm{c}w} = {A}_{\mathrm{o}w}+{B}_{\mathrm{o}w}K$, $K = [K_1 \ K_2]$,  ${B}_{\mathrm{c}w} = {N}_w+{B}_{\mathrm{o}w}K_3$, $ C_{\mathrm{c}w}=[{C}_w \ 0]+{D}_wK$, ${D}_{\mathrm{c}w} = {F}_w+{D}_wK_3$, 
\begin{align}\nonumber
    {A}_{\mathrm{o}w} &= \begin{bmatrix}
        {A}_w & 0 \\ G_2{C}_w & G_1
    \end{bmatrix}, \ {B}_{\mathrm{o}w} = \begin{bmatrix}
        {B}_w \\ G_2{D}_w
    \end{bmatrix}, \ 
    {N}_w =  \begin{bmatrix}
        {E}_w \\ G_2{F}_w
    \end{bmatrix}. \nonumber
\end{align}
The nominal (i.e., $w=0$) counterparts of the matrices $A_{\mathrm{c}w}$ and $B_{\mathrm{c}w}$ are denoted by $A_\mathrm{c}$ and $B_\mathrm{c}$.

We will refer to the following conditions while investigating the solvability of Problem \ref{prb:CRORP}.

\begin{condition}\label{ass:A0antiStable}
The matrix $A_0$ is diagonalizable, and all the eigenvalues of $A_0$ have zero real parts for continuous-time systems (respectively, modulus $1$ for discrete-time systems)\footnote{In Lemma \ref{lmm:keylemma}, $A_0$ can be allowed to be anti-stable; that is, all the eigenvalues of $A_0$ lie in the closed right half complex plane for continuous-time systems (respectively, have modulus greater than or equal to $1$ for discrete-time systems). However, since we require the exogenous signal $v$ not to leave the nondegenerate ellipsoidal set $\mathcal{V}$ throughout the letter, we discard the cases in which $v$ can be unbounded in Condition \ref{ass:A0antiStable}.}.  
\end{condition}

\begin{condition}\label{ass:pcopy}
The pair $(G_{1},G_{2})$ incorporates a $p$-copy internal model of $A_0$\footnote{The pair $(G_{1}, G_{2})$ takes the following form: $G_{1} = \mathrm{diag}(\alpha_{1},\ldots,\alpha_{p})$ and $G_{2} = \mathrm{diag}(\beta_{1},\ldots,\beta_{p})$, where for $\ell = 1, \ldots, p$,
$\alpha_{\ell}$ is a square matrix and $\beta_{\ell}$ is a column vector such that the minimal polynomial of $A_0$ is equal to
the characteristic polynomial of $\alpha_{\ell}$ and the pair $(\alpha_{\ell}, \beta_{\ell})$ is controllable for continuous-time systems and reachable for discrete-time systems.}.
\end{condition}


Conditions above are standard in linear output regulation theory (e.g., see Chapter 1 in \cite{huang2004nonlinear}). The following lemma, which is tailored from Lemma 1.27 in \cite{huang2004nonlinear}, is key for the solvability of Problem \ref{prb:CRORP} for both continuous-time and discrete-time systems.

\begin{lemma}\label{lmm:keylemma}
    Let Conditions \ref{ass:A0antiStable} and \ref{ass:pcopy} hold. If the nominal closed-loop system matrix $A_\mathrm{c}$ is stable,
    then for any $w \in \mathcal{W}$, where 
    $\mathcal{W}$ is a neighborhood of $w = 0$ in which the closed-loop system matrix $A_{\mathrm{c}w}$ is stable, there exists a unique matrix $X_{\mathrm{c}w}$ that satisfies the matrix equations  
    \begin{align}\label{eq:clsdloopreg} \nonumber
            X_{\mathrm{c}w}A_{0} &= {A}_{\mathrm{c}w}X_{\mathrm{c}w} + B_{\mathrm{c}w} \\
            0 &= {C}_{\mathrm{c}w}X_{\mathrm{c}w} + D_{\mathrm{c}w}.
        \end{align}
\end{lemma}
\begin{remark}\label{rmk:Existenceofw}
It is well known that if the nominal closed-loop system matrix $A_\mathrm{c}$ is stable, then there exists a neighborhood $\mathcal{W}$ of $w = 0$ such
that $A_{\mathrm{c}w}$ is stable for all $w \in \mathcal{W}$.
\end{remark}

Finally, we define the following ellipsoidal sets for the controller state $z$ and the closed-loop system state $\zeta$
\begin{align}\nonumber
    \mathcal{Q} &= \{z \in \mathbb{R}^{n_z} \ | \ z^\mathrm{T}Q z \leq 1 \}, \ Q \in \mathbb{S}^{n_z}_{++} \nonumber \\
    \Omega &= \{\zeta \in \mathbb{R}^{n+n_z} \ | \  \zeta^\mathrm{T} \Lambda \zeta \leq 1 \}, \ \Lambda = \mathrm{diag}(\Gamma,Q).  \nonumber
\end{align}

\section{Main Result for CRORP: Continuous Time}\label{sec:Continuous}

The following lemma employs Nagumo's invariance theorem to present a sufficient condition for the existence of an uncertainty set $\mathcal{W}$ containing $0$ in its interior, in which the set $\Omega$ is robustly positively invariant for the closed-loop system and the set $\mathcal{V}$. The sufficient condition is that the derivative of the quadratic function associated with the ellipsoidal set $\Omega$, i.e., $\zeta^\mathrm{T}\Lambda\zeta$, is negative along the trajectories of the nominal closed-loop system at every boundary point of $\Omega$ for all $v \in \mathcal{V}$.

\begin{lemma}\label{lmm:ct_strict_invariance}
    If $\zeta^\mathrm{T}\Lambda(A_\mathrm{c}\zeta + B_\mathrm{c}v) < 0$ for all $\zeta \in \partial\Omega$ and $v \in \mathcal{V}$, then there exists a neighborhood $\mathcal{W}$ of $w = 0$ such that,  for any $w \in \mathcal{W}$, $\Omega$ is robustly positively invariant for the closed-loop system \eqref{eq:augmented_clloop_system} and the set $\mathcal{V}$\footnote{Here, and likewise in Lemma \ref{lmm:ct_time_RPI}, $v$ is not necessarily generated from \eqref{eq:exosystem}.}.
\end{lemma}
\begin{IEEEproof}
    The tangent halfspace to the ellipsoidal set $\Omega$ at $\zeta \in \partial\Omega$ is $\mathcal{T}_\Omega(\zeta) = \{\eta \ | \ 2\zeta^\mathrm{T}\Lambda \eta \leq 0 \}$. By Nagumo's invariance theorem of interest (e.g., see Theorem 4.10 in \cite{Blanchinisbook},  Corollary 3 in \cite{pordal2025ellipsoidalsettheoreticdesignrobust}, and Section 4.1 in \cite{Blanchinisurvey}), given  $w \in \mathbb{R}^{\bar n}$, $\Omega$ is robustly positively invariant for the closed-loop system \eqref{eq:augmented_clloop_system}
    and the set $\mathcal{V}$ if
        \begin{align}\nonumber
        h(\zeta,v) \triangleq 2\zeta^\mathrm{T}\Lambda( A_{\mathrm{c}w} \zeta +  B_{\mathrm{c}w}v) \leq 0
    \end{align}
for all $\zeta \in \partial\Omega$ and $v \in \mathcal{V}$. Note that for any $w \in \mathbb{R}^{\bar n}$, $A_{\mathrm{c}w}$ and $B_{\mathrm{c}w}$ can be written as $A_{\mathrm{c}w} = A_\mathrm{c} + \Delta A_\mathrm{c}$ and  $B_{\mathrm{c}w} = B_\mathrm{c} + \Delta B_\mathrm{c}$, where 
    \begin{align}\nonumber
        \Delta A_\mathrm{c} &= \begin{bmatrix}
        \Delta{A} + \Delta{B}K_1 & \Delta{B}K_2 \\ G_2 (\Delta{C} + \Delta{DK_1}) & G_2 \Delta{DK_2}
    \end{bmatrix} \\
    \Delta B_\mathrm{c} &= \begin{bmatrix}
        \Delta{E} + \Delta{B}K_3 \\ G_2 (\Delta{F} + \Delta{D}K_3) 
    \end{bmatrix}. \nonumber
    \end{align}
Therefore,
$h$ can be equivalently written as 
   \begin{align}\label{eq:h_def}\nonumber
        h(\zeta,v) = \hspace{0.1cm} &\zeta^\mathrm{T}(A_\mathrm{c}^\mathrm{T}\Lambda + \Lambda A_\mathrm{c})\zeta + 2\zeta^\mathrm{T}\Lambda B_\mathrm{c}v \\ &+ 2\zeta^\mathrm{T}\Lambda( \Delta A_\mathrm{c} \zeta +  \Delta B_\mathrm{c}v).
    \end{align}
    
By the hypothesis,
    \begin{align} \nonumber
        \ell(\zeta,v) \triangleq \zeta^\mathrm{T}(A_\mathrm{c}^\mathrm{T}\Lambda + \Lambda A_\mathrm{c})\zeta + 2\zeta^\mathrm{T}\Lambda B_\mathrm{c}v < 0
    \end{align}
for all $\zeta \in \partial\Omega$ and $v \in \mathcal{V}$. 
 Since $\ell$ is a continuous function on the compact set $\partial\Omega \times \mathcal{V}$, it achieves its maximum over $\partial\Omega \times \mathcal{V}$. Accordingly, 
 \begin{align}\label{eq_upper_bounding_ell}
     \epsilon \triangleq  - \max_{(\zeta, v) \in \partial \Omega \times \mathcal{V}} \ell(\zeta, v) >0. 
 \end{align} 
 For any $w \in \mathbb{R}^{\bar n}$, we have 
\begin{align}\label{eq_uncertain_term_bounding}
& 2\zeta^\mathrm{T}\Lambda( \Delta A_\mathrm{c} \zeta +  \Delta B_\mathrm{c}v) \nonumber  \\
&  \leq  2 \|\Lambda\| \|\zeta\|^2 \|\Delta A_\mathrm{c}\| + 2 \|\Lambda\| \|\zeta\|  \|v\|  \|\Delta B_\mathrm{c}\|\nonumber  \\
& \leq \frac{2 \lambda_{\mathrm{max}}(\Lambda)}{\lambda_{\mathrm{min}}(\Lambda)}  \|\Delta A_\mathrm{c}\| +  \frac{2 \lambda_{\mathrm{max}}(\Lambda)}{\lambda_{\mathrm{min}}^{1/2}(\Lambda) \lambda_{\mathrm{min}}^{1/2}(R)  } \|\Delta B_\mathrm{c}\| 
\end{align}
for all $\zeta \in \partial\Omega$ and $v \in \mathcal{V}$, where the last inequality follows from the facts $\|\Lambda\| = \lambda_{\mathrm{max}}(\Lambda)$,  $\lambda_{\mathrm{min}}(\Lambda) \|\zeta\|^2 \leq 1$ for all $\zeta \in \partial\Omega$, and $\lambda_{\mathrm{min}}(R) \| v \|^2 \leq 1$  for all $v \in \mathcal{V}$. Since $\Delta A_{\mathrm{c}}$ and $\Delta B_{\mathrm{c}}$ are continuous in $w$ and vanish at $w=0$, there exists a neighborhood $\mathcal{W}$ of $w=0$ such that $\|\Delta A_{\mathrm{c}}\| < \delta$ and $\|\Delta B_{\mathrm{c}}\| < \delta$ for all $w \in \mathcal{W}$, where 
\begin{align}\nonumber
\delta =  \min \left( \frac{\lambda_{\mathrm{min}}(\Lambda) \epsilon  }{4\lambda_{\mathrm{max}}(\Lambda)} 
,   \frac{\lambda_{\mathrm{min}}^{1/2} (\Lambda) \lambda_{\mathrm{min}}^{1/2}(R) \epsilon}{4\lambda_{\mathrm{max}}(\Lambda)}    \right) >0. 
\end{align}
Using this and \eqref{eq_uncertain_term_bounding}, one can show that, for any $w \in \mathcal{W}$, 
\begin{align}\label{eq:uncertain+term_bounding_final}
    & 2\zeta^\mathrm{T}\Lambda( \Delta A_\mathrm{c} \zeta +  \Delta B_\mathrm{c}v) < \epsilon  
\end{align}
for all $\zeta \in \partial\Omega$ and $v \in \mathcal{V}$. We now conclude from \eqref{eq:h_def}, \eqref{eq_upper_bounding_ell}, and \eqref{eq:uncertain+term_bounding_final} that, for any $w\in \mathcal{W}$, $h(\zeta, v) \leq 0$ for all $\zeta \in \partial\Omega$ and $v \in \mathcal{V}$.
\end{IEEEproof}

Building on Lemma \ref{lmm:ct_strict_invariance}, the following result provides an LMI condition for the existence of an uncertainty set $\mathcal{W}$ containing $0$ in its interior, in which the closed-loop system matrix is Hurwitz and the set $\Omega$ is robustly positively invariant for the closed-loop system and the set $\mathcal{V}$.
\begin{lemma}\label{lmm:ct_time_RPI}
    If there exist $K$, $K_3$, and $\alpha >0$ such that
    \begin{align}\label{eq:ct_solvability_LMI}
       \begin{bmatrix}
            A_\mathrm{c}^\mathrm{T}\Lambda + \Lambda A_\mathrm{c}+\alpha \Lambda & \Lambda B_\mathrm{c} \\ B_\mathrm{c}^\mathrm{T}\Lambda & -\alpha R
       \end{bmatrix} \prec 0,
    \end{align}
    then 
    there exists a neighborhood $\mathcal{W}$ of $w = 0$ such that, for any $w \in \mathcal{W}$,
    \begin{enumerate}
        \item[(i)] the closed-loop system matrix $A_{\mathrm{c}w}$ is Hurwitz and
        \item[(ii)] the set $\Omega$ is robustly positively invariant for the closed-loop system \eqref{eq:augmented_clloop_system} and the set $\mathcal{V}$.
    \end{enumerate}
\end{lemma}
\begin{IEEEproof}
Since $\alpha > 0 $ and $\Lambda \succ 0$, we have   $ A_\mathrm{c}^\mathrm{T}\Lambda + \Lambda A_\mathrm{c} \prec 0$ necessarily by \eqref{eq:ct_solvability_LMI}. Hence,  $A_\mathrm{c}$ is Hurwitz. By Remark \ref{rmk:Existenceofw}, there exists a neighborhood $\mathcal{W}_\mathrm{1}$ of $w = 0$ such that $A_{\mathrm{c}w}$ is Hurwitz for all $w \in \mathcal{W}_\mathrm{1}$.

Note that, for any $\zeta \in \mathbb{R}^{n+n_z}$ and $v\in \mathbb{R}^{q}$,
\begin{align}\label{eq_lyap_derivative_lem3}
    2\zeta^\mathrm{T}\Lambda(A_\mathrm{c}\zeta + B_\mathrm{c}v) = \zeta^\mathrm{T}( A_\mathrm{c}^\mathrm{T}\Lambda + \Lambda A_\mathrm{c})\zeta + 2\zeta^\mathrm{T}\Lambda B_\mathrm{c}v.
\end{align}
 Using $(a-b)^\mathrm{T}(a-b) \geq 0$ with $a = \alpha^{-1/2}R^{-1/2}B_\mathrm{c}^\mathrm{T}\Lambda\zeta$ and $b = \alpha^{1/2} R^{1/2}v$, one can also establish that
    \begin{align}\label{eq:cross_term_bounding}
     2\zeta^\mathrm{T}\Lambda B_\mathrm{c}v \leq \frac{1}{\alpha} \zeta^\mathrm{T}\Lambda B_\mathrm{c}R^{-1}B_\mathrm{c}^\mathrm{T}\Lambda \zeta + \alpha v^\mathrm{T}Rv.   
    \end{align}
    Since  $ \zeta^\mathrm{T}\Lambda \zeta = 1$ for all $\zeta \in \partial\Omega$ and $ v^\mathrm{T}Rv \leq 1$ for all $v \in \mathcal{V}$, $\alpha v^\mathrm{T}Rv \leq \alpha \zeta^\mathrm{T}\Lambda \zeta$ for all $\zeta \in \partial\Omega$ and $v \in \mathcal{V}$.  In conjunction with \eqref{eq_lyap_derivative_lem3} and \eqref{eq:cross_term_bounding}, this implies that 
    \begin{align}\nonumber
 & \zeta^\mathrm{T}\Lambda(A_\mathrm{c}\zeta + B_\mathrm{c}v) \\& \leq \frac{1}{2}\zeta^\mathrm{T}(\underbrace{A_\mathrm{c}^\mathrm{T}\Lambda + \Lambda A_\mathrm{c} + \frac{1}{\alpha} \Lambda B_\mathrm{c} R^{-1} B_\mathrm{c}^\mathrm{T}\Lambda + \alpha\Lambda}_\text{\small $\triangleq -S$})\zeta \label{eq_S_involved_ineq}
\end{align}
    for all $\zeta \in \partial\Omega$ and $v \in \mathcal{V}$. By the LMI \eqref{eq:ct_solvability_LMI} and the Schur complement of $-\alpha R$ in the matrix in \eqref{eq:ct_solvability_LMI}, $S \succ 0$. Thus, we now infer from \eqref{eq_S_involved_ineq} that $\zeta^\mathrm{T}\Lambda(A_\mathrm{c}\zeta + B_\mathrm{c}v) < 0 $ for all $\zeta \in \partial\Omega$ and $v \in \mathcal{V}$. 
    By Lemma \ref{lmm:ct_strict_invariance},  there exists a neighborhood $\mathcal{W}_2$ of $w=0$ such that, for any $w\in \mathcal{W}_2$, $\Omega$ is robustly positively invariant for the closed-loop system \eqref{eq:augmented_clloop_system} and the set $\mathcal{V}$.
    
Since the sets $\mathcal{W}_1$ and $\mathcal{W}_2$ are both neighborhoods of $w = 0$, 
their intersection has a nonempty interior that contains $w=0$. Therefore, there exists a neighborhood $\mathcal{W}$ of $w = 0$ such that for any $w \in \mathcal{W}$, both (i) and (ii) in Lemma \ref{lmm:ct_time_RPI} hold. 
\end{IEEEproof}

Finally, Lemmas \ref{lmm:keylemma} and \ref{lmm:ct_time_RPI} pave the way for the main result of the continuous-time CRORP.

\begin{theorem}\label{thm:main_thm_ct_time}
    Let Conditions \ref{ass:A0antiStable} and \ref{ass:pcopy} hold. If there exist $K$, $K_3$, $\alpha>0$, and $\beta \in [0,1]$ such that the LMIs \eqref{eq:ct_solvability_LMI} and
    \begin{align}\label{eq:ctrl_related_LMIs}
        \begin{bmatrix}
           \beta \Lambda & 0 & K^\mathrm{T} \\ 0 & (1-\beta) R & K_3^\mathrm{T} \\
            K & K_3 & \Pi^{-1}
        \end{bmatrix} \succeq 0
    \end{align}
    hold, then 
    the dynamic state-feedback control law \eqref{eq:controllaw} solves Problem \ref{prb:CRORP}. 
\end{theorem}

\begin{IEEEproof}
 By Lemma \ref{lmm:ct_time_RPI}, there exists a neighborhood $\mathcal{W}$ of $w = 0$ such that, for any $w \in \mathcal{W}$, both (i) and (ii) in Lemma \ref{lmm:ct_time_RPI}  hold. Let $w \in \mathcal{W} $ be given. Also, fix $x(0) \in \mathcal{X}$ and $z(0) = 0$, and let the exogenous signal generated by \eqref{eq:exosystem} satisfy   $v(t) \in \mathcal{V}$ for all $t \geq 0$. 
 Since $x(0) \in \mathcal{X}$ and $z(0) = 0$, 
 $\zeta^{\mathrm{T}}(0)\Lambda\zeta(0) =  x^\mathrm{T}(0)\Gamma x(0) \leq 1$, and hence, $ \zeta(0) \in \Omega$. We now deduce from the set $\Omega$ being robustly positively invariant for the closed-loop system \eqref{eq:augmented_clloop_system} and the set $\mathcal{V}$ that $\zeta(t) \in \Omega$ for all $t > 0$. That is, $\zeta^\mathrm{T}(t)\Lambda \zeta(t) = x^\mathrm{T}(t)\Gamma x(t) + z^\mathrm{T}(t) Q z(t) \leq 1$ for all $t >0 $, implying $x^\mathrm{T}(t)\Gamma x(t) \leq 1$ for all $t> 0$, and hence, $x(t) \in \mathcal{X}$ for all $t > 0$.

 Next, the LMI \eqref{eq:ctrl_related_LMIs} implies
     \begin{equation}\label{eq:input_constraint_LMI_implication}
         \begin{bmatrix}
         \beta  \Lambda & 0 \\ 0 & (1-\beta)R
         \end{bmatrix} \succeq \begin{bmatrix}
             K^\mathrm{T} \\ K_3^\mathrm{T}
         \end{bmatrix} \Pi \begin{bmatrix}
             K & K_3
         \end{bmatrix}
     \end{equation}
     by the Schur complement of $ \Pi^{-1}$ in the matrix in \eqref{eq:ctrl_related_LMIs}. 
     Note that $u = K \zeta + K_3v$. Therefore,
     \begin{align}\nonumber 
       u^\mathrm{T} \Pi u =  \begin{bmatrix}
           \zeta^{\mathrm{T}} & v^{\mathrm{T}} 
       \end{bmatrix} \begin{bmatrix}
             K^\mathrm{T} \\ K_3^\mathrm{T}
         \end{bmatrix} \Pi \begin{bmatrix}
             K & K_3
         \end{bmatrix} \begin{bmatrix}
             \zeta \\ v
         \end{bmatrix}. 
     \end{align}
  This, together with  \eqref{eq:input_constraint_LMI_implication}, implies that 
  \begin{align}\nonumber
   u^\mathrm{T}(t) \Pi u(t) \leq \beta  \zeta^{\mathrm{T}}(t) \Lambda \zeta(t) + (1-\beta) v^{\mathrm{T}}(t)R  v(t) 
  \end{align}
for all $t\geq 0$. Since $0 \leq \beta \leq 1$ and  $\zeta(t) \in \Omega$ and $v(t) \in \mathcal{V}$ for all $t \geq 0$, we infer that $u^\mathrm{T}(t) \Pi u(t) \leq 1$ for all $t \geq 0$, and hence, $u(t) \in \mathcal{U}$ for all $t \geq 0$. 
  
 Under  Conditions \ref{ass:A0antiStable} and \ref{ass:pcopy}, there exists a unique solution $X_{\mathrm{c}w}$ to the matrix equations \eqref{eq:clsdloopreg} as per Lemma \ref{lmm:keylemma}. Let $\Bar{\zeta} = \zeta - X_{\mathrm{c}w}v$. Then, using \eqref{eq:exosystem}, \eqref{eq:augmented_clloop_system}, and \eqref{eq:clsdloopreg}, one can show that
     \begin{align}
         \rho{\Bar{\zeta}} = A_{\mathrm{c}w}\Bar{\zeta} \ \text{ and } \ e = C_{\mathrm{c}w}\Bar{\zeta}.
     \end{align}  
In conjunction with the Hurwitzness of $A_{\mathrm{c}w}$, this implies that 
$\lim_{t\to\infty} \bar \zeta(t) = 0$, and hence, $\lim_{t\to\infty} e(t) = 0$. 
\end{IEEEproof}

\section{Main Result for CRORP: Discrete Time}


This section provides discrete-time counterparts of the results in Section \ref{sec:Continuous} within a similar logical flow. Specifically, Lemma \ref{lmm:dt_strict_invariance}, Lemma \ref{lmm:dt_time_RPI}, and Theorem \ref{thm:main_thm_dt_time} are discrete-time counterparts of Lemma \ref{lmm:ct_strict_invariance}, Lemma \ref{lmm:ct_time_RPI}, and Theorem \ref{thm:main_thm_ct_time}, respectively. The sufficient condition in Lemma \ref{lmm:dt_strict_invariance} for the existence of an uncertainty set of interest can be interpreted as: for any $\zeta \in \Omega$, the successor state of the nominal closed-loop system is in the interior of $\Omega$ for all $v \in \mathcal V$.







\begin{lemma}\label{lmm:dt_strict_invariance}
    If $A_\mathrm{c}\zeta + B_\mathrm{c}v \in \mathrm{int} \Omega$ for all $\zeta \in \Omega$ and $v \in \mathcal{V}$, then there exists a neighborhood $\mathcal{W}$ of $w = 0$ such that,  for any $w \in \mathcal{W}$, $\Omega$ is robustly positively invariant for the closed-loop system \eqref{eq:augmented_clloop_system} and the set $\mathcal{V}$.
\end{lemma}
\begin{IEEEproof}
It can be easily derived from Definition \ref{def:RPIset} that,  given  $w \in \mathbb{R}^{\bar n}$, $\Omega$ is robustly positively invariant for the closed-loop system \eqref{eq:augmented_clloop_system}
    and the set $\mathcal{V}$ if $A_{\mathrm{c}w}\zeta + B_{\mathrm{c}w}v \in \Omega$ for all $\zeta \in \Omega$ and $v \in \mathcal{V}$ (e.g., see Theorem 4 in  \cite{8945668}), or equivalently,  
    \begin{align}
        h(\zeta, v) \triangleq      (A_{\mathrm{c}w}\zeta + B_{\mathrm{c}w}v)^\mathrm{T} \Lambda (A_{\mathrm{c}w}\zeta + B_{\mathrm{c}w}v) \leq 1 
    \end{align}
for all    $\zeta \in \Omega$ and $v \in \mathcal{V}$. As in the proof of Lemma \ref{lmm:ct_strict_invariance}, using   $A_{\mathrm{c}w} = A_\mathrm{c} + \Delta A_\mathrm{c}$ and  $B_{\mathrm{c}w} = B_\mathrm{c} + \Delta B_\mathrm{c}$ , $h$ can be equivalently written as 
\begin{align}\label{eq_h_rewritten_dt_lem4}
    h(\zeta, v) =& (A_\mathrm{c}\zeta + B_\mathrm{c} v)^\mathrm{T} \Lambda (A_\mathrm{c}\zeta + B_\mathrm{c} v) \nonumber 
     \\ & + 
    2(A_\mathrm{c}\zeta + B_\mathrm{c} v)^\mathrm{T} \Lambda (\Delta A_\mathrm{c}\zeta + \Delta B_\mathrm{c} v) \nonumber  
    \\ & + (\Delta A_\mathrm{c}\zeta + \Delta B_\mathrm{c} v)^\mathrm{T} \Lambda (\Delta A_\mathrm{c}\zeta + \Delta B_\mathrm{c} v). 
\end{align}

By the hypothesis, 
\begin{align}\nonumber
    \ell(\zeta, v) \triangleq (A_\mathrm{c}\zeta + B_\mathrm{c} v)^\mathrm{T} \Lambda (A_\mathrm{c}\zeta + B_\mathrm{c} v) < 1
\end{align}
for all $\zeta \in \Omega$ and $v \in \mathcal{V}$. Since $g(\zeta, v) \triangleq 1 - \ell(\zeta, v)$ is a continuous function on the compact set $\Omega \times \mathcal{V}$, it attains its minimum on $\Omega \times \mathcal{V}$. 
Accordingly, 
 \begin{align}\label{eq_lower_bounding_g}
     \epsilon \triangleq   \min_{(\zeta, v) \in \Omega \times \mathcal{V}} g(\zeta, v) >0. 
 \end{align} 
Note that $\|\Lambda\| = \lambda_{\mathrm{max}}(\Lambda)$,  $\lambda_{\mathrm{min}}(\Lambda) \|\zeta\|^2 \leq 1$ for all $\zeta \in \Omega$, and $\lambda_{\mathrm{min}}(R) \| v \|^2 \leq 1$  for all $v \in \mathcal{V}$. 
For any
$w \in \mathbb{R}^{\bar n}$, using also the Cauchy-Schwarz inequality and  norm properties, 
one can show that
\begin{align}\label{eq_mix_term_upper_bounding_lem4}
& 2(A_\mathrm{c}\zeta + B_\mathrm{c} v)^\mathrm{T} \Lambda (\Delta A_\mathrm{c}\zeta + \Delta B_\mathrm{c} v) \nonumber \\ & \leq a \| \Delta A_{\mathrm{c}} \|  +  b \| \Delta B_{\mathrm{c}} \|   \end{align}
for all  $\zeta \in \Omega$ and $v \in \mathcal{V}$, where 
\begin{align}
   a &=  2 \lambda_{\mathrm{max}}(\Lambda) \left( \frac{\| {A_{\mathrm{c}}} \|}{\lambda_{\mathrm{min}}(\Lambda)} +\frac{\| B_{\mathrm{c}}\|}{\lambda_{\mathrm{min}}^{1/2}(\Lambda) \lambda_{\mathrm{min}}^{1/2}(R)} \right)  \nonumber \\
   b &= 2 \lambda_{\mathrm{max}}(\Lambda) \left( \frac{\| A_{\mathrm{c}}\|}{\lambda_{\mathrm{min}}^{1/2}(\Lambda) \lambda_{\mathrm{min}}^{1/2}(R)} +  \frac{\| {B_{\mathrm{c}}} \|}{\lambda_{\mathrm{min}}(R)}  \right).   \nonumber   \nonumber 
\end{align}
 Similarly, for any
$w \in \mathbb{R}^{\bar n}$, we have
 \begin{align}\label{eq_uncertain_term_bounding_square_lem4}
  &   (\Delta A_\mathrm{c}\zeta + \Delta B_\mathrm{c} v)^\mathrm{T} \Lambda (\Delta A_\mathrm{c}\zeta + \Delta B_\mathrm{c} v)  \nonumber \\ & \leq \left( \frac{\lambda_{\mathrm{max}}^{1/2}(\Lambda) }{\lambda_{\mathrm{min}}^{1/2}(\Lambda)} \| \Delta A_{\mathrm{c}} \|  + \frac{\lambda_{\mathrm{max}}^{1/2}(\Lambda) }{\lambda_{\mathrm{min}}^{1/2}(R)}  \| \Delta B_{\mathrm{c}} \| \right)^2  
 \end{align}
 for all  $\zeta \in \Omega$ and $v \in \mathcal{V}$. Since $\Delta A_{\mathrm{c}}$ and $\Delta B_{\mathrm{c}}$ are continuous in $w$ and vanish at $w=0$, there exists a neighborhood $\mathcal{W}$ of $w=0$ such that $\|\Delta A_{\mathrm{c}}\| < \delta$ and $\|\Delta B_{\mathrm{c}}\| < \delta$ for all $w \in \mathcal{W}$, where 
\begin{equation}\nonumber
\delta \hspace{-0.04cm} =  \hspace{-0.05cm} \begin{cases}
 \min  \left( \frac{\epsilon^{1/2}\lambda_{\mathrm{min}}^{1/2}(\Lambda)}{2\lambda_{\mathrm{max}}^{1/2}(\Lambda)}  ,  \frac{\epsilon^{1/2}\lambda_{\mathrm{min}}^{1/2}(R)}{2\lambda_{\mathrm{max}}^{1/2}(\Lambda)}  \right)  \  \hspace{-0.05cm }\mathrm{if}\  \| A_{\mathrm{c}}\| =  \| B_{\mathrm{c}}\| =0 \\ 
\min  \left( \frac{\epsilon}{3a}, \frac{\epsilon}{3b}, \frac{\epsilon^{1/2}\lambda_{\mathrm{min}}^{1/2}(\Lambda)}{12^{1/2}\lambda_{\mathrm{max}}^{1/2}(\Lambda)}  ,  \frac{\epsilon^{1/2}\lambda_{\mathrm{min}}^{1/2}(R)}{12^{1/2}\lambda_{\mathrm{max}}^{1/2}(\Lambda)}  \right)    \  \hspace{-0.05cm } \mathrm{otherwise}.
\end{cases}
\end{equation} 
Using this, \eqref{eq_h_rewritten_dt_lem4}, the definition of the function $\ell$, \eqref{eq_mix_term_upper_bounding_lem4}, and \eqref{eq_uncertain_term_bounding_square_lem4}, one can show that, for any $w \in \mathcal{W}$,
\begin{align}\label{eq_h_subtract_l_lemm4}
    h(\zeta, v) - \ell(\zeta,v) < \epsilon 
\end{align}
for all $\zeta \in \Omega$ and $v \in \mathcal{V}$. We can now conclude from  $\ell(\zeta,v) = 1 - g(\zeta, v)$, \eqref{eq_lower_bounding_g}, and \eqref{eq_h_subtract_l_lemm4} that, for any $w \in \mathcal{W}$,  $h(\zeta, v) \leq 1$  for all $\zeta \in \Omega$ and $v \in \mathcal{V}$.
 \end{IEEEproof}

\begin{lemma}\label{lmm:dt_time_RPI}
    If there exist $K$, $K_3$, and $\gamma \in (0,1)$ such that
    \begin{align}\label{eq:dt_solvability_LMI}
       \begin{bmatrix}
            \gamma\Lambda & 0 & A_\mathrm{c}^\mathrm{T} \\ 
            0 & (1-\gamma)R & B_\mathrm{c}^\mathrm{T} \\ 
            A_\mathrm{c} & B_\mathrm{c} & \Lambda^{-1}
       \end{bmatrix} \succ 0,
    \end{align}
    then 
    there exists a neighborhood $\mathcal{W}$ of $w = 0$ such that, for any $w \in \mathcal{W}$,
    \begin{enumerate}
        \item[(i)] the closed-loop system matrix $A_{\mathrm{c}w}$ is Schur and
        \item[(ii)] the set $\Omega$ is robustly positively invariant for the closed-loop system \eqref{eq:augmented_clloop_system} and the set $\mathcal{V}$.
    \end{enumerate}
\end{lemma}
\begin{IEEEproof}
    By the Schur complement of $ \Lambda^{-1}$ in the matrix in \eqref{eq:dt_solvability_LMI}, the LMI \eqref{eq:dt_solvability_LMI} implies 
\begin{align}\label{eq:dt_solvability_LMI_equiv1}
       \begin{bmatrix}
            A_\mathrm{c}^\mathrm{T}\Lambda A_\mathrm{c} & A_\mathrm{c}^\mathrm{T}\Lambda B_\mathrm{c} \\ B_\mathrm{c}^\mathrm{T}\Lambda A_\mathrm{c} & B_\mathrm{c}^\mathrm{T}\Lambda B_\mathrm{c}
       \end{bmatrix} \prec \begin{bmatrix}
           \gamma\Lambda & 0 \\ 0 & (1-\gamma)R
       \end{bmatrix}.
    \end{align}
     Since $0<\gamma < 1$ and $\Lambda \succ 0$, we have  $A_\mathrm{c}^\mathrm{T}\Lambda A_\mathrm{c} - \Lambda \prec 0$ necessarily by \eqref{eq:dt_solvability_LMI_equiv1}. Hence,  $A_\mathrm{c}$ is Schur. By Remark \ref{rmk:Existenceofw}, there exists a neighborhood $\mathcal{W}_\mathrm{1}$ of $w = 0$ such that $A_{\mathrm{c}w}$ is Schur for all $w \in \mathcal{W}_\mathrm{1}$.

 Note that, for any $\zeta \in \mathbb{R}^{n+n_z}$ and $v\in \mathbb{R}^{q}$,    
\begin{align}\label{eq_discrete_lyap_related_lem5}
    &  (A_{\mathrm{c}}\zeta + B_{\mathrm{c}}v)^\mathrm{T} \Lambda (A_{\mathrm{c}}\zeta + B_{\mathrm{c}}v) \nonumber \\ 
    & =\zeta^{\mathrm{T}}A_{\mathrm{c}}^{\mathrm{T}}\Lambda A_{\mathrm{c}} \zeta  + 2 \zeta^{\mathrm{T}}A_{\mathrm{c}}^{\mathrm{T}}\Lambda B_{\mathrm{c}} v + v^{\mathrm{T}} B_{\mathrm{c}}^{\mathrm{T}}\Lambda  B_{\mathrm{c}}v. 
\end{align}
By \eqref{eq:dt_solvability_LMI_equiv1}, we also have 
\begin{align}\label{eq:dt_LMI_quadratic}\nonumber
&\zeta^\mathrm{T}A_\mathrm{c}^\mathrm{T}\Lambda A_\mathrm{c} \zeta + 2\zeta^\mathrm{T}A_\mathrm{c}^\mathrm{T}\Lambda B_\mathrm{c}v + v^\mathrm{T}B_\mathrm{c}^\mathrm{T}\Lambda B_\mathrm{c}v \\ 
    &< \gamma \zeta^\mathrm{T}\Lambda \zeta + (1-\gamma)v^\mathrm{T}R v. 
    \end{align}
   Since $0<\gamma<1$,   $ \zeta^\mathrm{T}\Lambda \zeta \leq  1$ for all $\zeta \in \Omega$, and $ v^\mathrm{T}Rv \leq 1$ for all $v \in \mathcal{V}$, the following inequality
   \begin{align}\nonumber
       \gamma \zeta^\mathrm{T}\Lambda \zeta + (1-\gamma)v^\mathrm{T}R v \leq  1 
   \end{align}
 holds  for all  $\zeta \in \Omega$ and $v \in \mathcal{V}$. This, together with  \eqref{eq_discrete_lyap_related_lem5} and \eqref{eq:dt_LMI_quadratic}, implies 
 that 
 \begin{align}\nonumber
    (A_\mathrm{c}\zeta + B_\mathrm{c}v)^\mathrm{T}\Lambda (A_\mathrm{c}\zeta+B_\mathrm{c})v <  1   
    \end{align}
    for all  $\zeta \in \Omega$ and $v \in \mathcal{V}$. Thus, $A_\mathrm{c}\zeta + B_\mathrm{c}v \in \mathrm{int} \Omega$ for all $\zeta \in \Omega$ and $v \in \mathcal{V}$.
    By Lemma \ref{lmm:dt_strict_invariance}, there exists a neighborhood $\mathcal{W}_2$ of $w = 0$ such that, for any $w \in \mathcal{W}_2$, the set $\Omega$ is robust positively invariant for the closed-loop system \eqref{eq:augmented_clloop_system} and the set $\mathcal{V}$. The existence of a neighborhood $\mathcal{W}$ of $w = 0$ such that, for any $w \in \mathcal{W}$, both (i) and (ii) in Lemma \ref{lmm:dt_time_RPI}  hold, follows from the same argument in the proof of Lemma \ref{lmm:ct_time_RPI}. 
\end{IEEEproof}

\begin{theorem}\label{thm:main_thm_dt_time}
    Let Conditions \ref{ass:A0antiStable} and \ref{ass:pcopy} hold. If there exist $K$, $K_3$, $\beta \in [0,1]$, and $\gamma \in (0,1)$ such that the LMIs \eqref{eq:ctrl_related_LMIs} and \eqref{eq:dt_solvability_LMI} hold, then
    the dynamic state-feedback control law \eqref{eq:controllaw} solves Problem \ref{prb:CRORP}. 
\end{theorem}
\begin{IEEEproof}
    Utilizing Lemma \ref{lmm:dt_time_RPI} instead of Lemma \ref{lmm:ct_time_RPI}, one can conclude the proof by following the steps in the proof of Theorem \ref{thm:main_thm_ct_time}.
\end{IEEEproof}
\begin{figure*}[!hb]
  \begin{center}
    \includegraphics[width=0.975\textwidth,trim=145 5 125 20,clip]{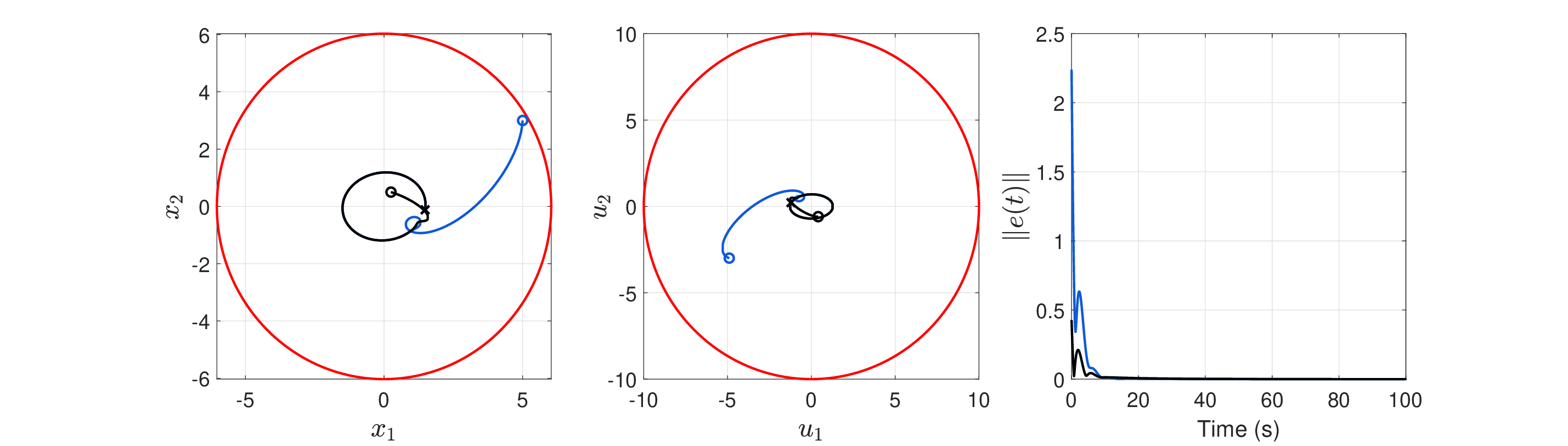} 
    \caption{State trajectories and control inputs in the respective phase planes, and tracking error norms over time, for two different initial conditions. Signals associated with each initial condition are shown with the same colors across all three plots. Boundaries of the respective constraint sets are shown in red, while ``o" and ``$\times$" marks  the initial and terminal points, respectively.}
  \label{fig:numex}
  \end{center}
\end{figure*}

\section{Illustrative Numerical Example}

We consider a continuous-time system in the form \eqref{eq:uncertain_system} with 
\begin{align}\nonumber
    A &= \begin{bmatrix}
          0.75 & 0.3 \\
              0 & 0.5
    \end{bmatrix}, \ B = \begin{bmatrix}
        1 & 0 \\ 0.2 & 1
    \end{bmatrix}, \ C = \begin{bmatrix}
        1 & 0 \\ 0 & 1
    \end{bmatrix}, 
\end{align}
$D = 0.5I_2$, $E = 0$, and $F = -I_2$, and the exosystem \eqref{eq:exosystem} with
\begin{equation}\nonumber
    A_0 = \begin{bmatrix}
        0 & 2\pi f \\ -2\pi f & 0
    \end{bmatrix},
\end{equation}
where $f = 0.01 \ \mathrm{Hz} $. Matrices defining the ellipsoidal constraint sets are $\Gamma = 0.0275I_2$, $\Pi = 0.01I_2$, $R = 0.75I_2$, and
\begin{equation}\nonumber
    Q = \begin{bmatrix}
   0.2729  & 0.1083 & 0.0108 & 0.0044 \\
   0.1083  & 0.1005 & 0.0044 & 0.0037 \\
   0.0108  & 0.0044 & 0.2561 & 0.1013 \\
   0.0044  & 0.0037 & 0.1013 & 0.0948
    \end{bmatrix}.   
\end{equation}

We let $G_1 = \mathrm{diag}(A_0,A_0)$ and $G_2 = \mathrm{diag}([0,1]^\mathrm{T},[0,1]^\mathrm{T})$. Note that Conditions \ref{ass:A0antiStable}--\ref{ass:pcopy} hold. To synthesize the control gains, we use CVX \cite{grant2014cvx} to model the semidefinite feasibility program yielded by the LMIs in Theorem \ref{thm:main_thm_ct_time} and we employ SDPT3 \cite{Toh1999SDPT3} to solve it in MATLAB.
The resulting control gains are 
\begin{align}\nonumber
    K_1 &= \begin{bmatrix}
         -1.0925 & -0.0467 \\
   0.0164 & -0.9913
    \end{bmatrix}, \ K_3 = \begin{bmatrix}
           0.7027 &  0.0016 \\
  -0.1014 & 0.7648
    \end{bmatrix}, \\ K_2 &= \begin{bmatrix}
        -1.2685 & -1.2659 & -0.0772 & -0.0805 \\
   0.1088 & 0.093 & -1.297 & -1.2904
    \end{bmatrix}.\nonumber
\end{align}

Uncertainties considered are given by $\Delta A = \Delta B = \Delta C = \Delta D = \Delta E = \Delta F = 0.1I_2$. Fig. \ref{fig:numex} includes two simulation results. First one is initiated with $v(0) = 1$, $z(0) = 0$, and $x^\mathrm{T}(0) = [5,3]$ and the second one is with $v(0) = 1$, $z(0) = 0$, and $x^\mathrm{T}(0) = [0.25,0.5]$. As promised by Theorem \ref{thm:main_thm_ct_time}, the state trajectories and the corresponding control inputs respect their constraints while the present disturbance is rejected and the reference is tracked by the output of the system despite the uncertainties, that is, the CRORP is solved. 

\section*{Conclusion}
Robust output regulation of both continuous-time and discrete-time uncertain LTI systems subject to ellipsoidal state and input constraints is considered. LMI conditions that lead to convex semidefinite feasibility programs for control parameter synthesis are provided. Extending the results of this letter to data-driven CRORP \cite{11664319} and robust cooperative output regulation of multi-agent systems \cite{SarsilmazIJC,koru2020cooperative,DTRCORP} subject to ellipsoidal constraints are potential directions for future research.


\section*{Acknowledgment}
This work was supported  in part through the Utah NASA Space Grant Consortium, NASA Grant \#80NSSC25M7088. The authors also thank Akalu Desta Teklu for his participation in initial discussions of this letter.

\bibliographystyle{IEEEtran}
\bibliography{Burak.bib}

@book{Blanchinisbook,
author="Blanchini, Franco
and Miani, Stefano",
title="Set-Theoretic Methods in Control",
year="2008",
publisher="Birkh{\"a}user Boston",
address="Boston, MA",
}

@article{Blanchinisurvey,
title = {Set invariance in control},
journal = {Automatica},
volume = {35},
number = {11},
pages = {1747-1767},
year = {1999},
issn = {0005-1098},
author = {F. Blanchini},
}

@article{VASSILAKI01061988,
author = {M. Vassilaki and J. C. Hennet and G. Bitsoris},
title = {Feedback control of linear discrete-time systems under state and control constraints},
journal = {Int. J. Control.},
volume = {47},
number = {6},
pages = {1727--1735},
year = {1988},
publisher = {Taylor \& Francis},
}

@ARTICLE{1104270,
  author={Gutman, P.-O. and Cwikel, M.},
  journal={IEEE Trans. Autom. Control.}, 
  title={Admissible sets and feedback control for discrete-time linear dynamical systems with bounded controls and states}, 
  year={1986},
  volume={31},
  number={4},
  pages={373-376},}

@article{MILANI19951489,
title = {Robust linear regulator design for discrete-time systems under polyhedral constraints},
journal = {Automatica},
volume = {31},
number = {10},
pages = {1489-1493},
year = {1995},
issn = {0005-1098},
author = {Basílio E.A. Milani and Andréa N. Carvalho}
}

@article{Kurzhanski_ellipsoidal,
author = { Kurzhanski, A. B. and  Varaiya, P.},
title = {Ellipsoidal Techniques for Reachability Under State Constraints},
journal = {SIAM J. Control Optim.},
volume = {45},
number = {4},
pages = {1369-1394},
year = {2006},
}

@misc{pordal2025ellipsoidalsettheoreticdesignrobust,
      title={Ellipsoidal Set-Theoretic Design of Robust Safety Filters for Constrained Linear Systems}, 
      author={Reza Pordal and Alireza Sharifi and Ali Baniasad},
      year={2025},
      eprint={2510.22790v2},
      archivePrefix={arXiv},
      primaryClass={eess.SY},
      url={https://arxiv.org/abs/2510.22790v2}, 
}

@article{KOTHARE19961361,
title = {Robust constrained model predictive control using linear matrix inequalities},
journal = {Automatica},
volume = {32},
number = {10},
pages = {1361-1379},
year = {1996},
issn = {0005-1098},
author = {Mayuresh V. Kothare and Venkataramanan Balakrishnan and Manfred Morari},
}

@article{BROOMS2001157,
title = {Constrained {MPC} for uncertain linear systems with ellipsoidal target sets},
journal = {Syst. Control Lett.},
volume = {44},
number = {3},
pages = {157-166},
year = {2001},
issn = {0167-6911},
author = {A.C. Brooms and B. Kouvaritakis and Y.I. Lee},
}

@book{IsidoriORP,
 author = {Byrnes, Christopher I.
and Priscoli, Francesco Delli
and Isidori, Alberto},
 publisher = {Birkh{\"a}user Boston},
 title = {Output Regulation of Nonlinear Systems},
 year = {1997}
}

@ARTICLE{1101137,
  author={Davison, E.},
  journal={IEEE Trans. Autom. Control.}, 
  title={The robust control of a servomechanism problem for linear time-invariant multivariable systems}, 
  year={1976},
  volume={21},
  number={1},
  pages={25-34}}

@ARTICLE{10502160,
  author={Zhu, Lijun and Chen, Zhiyong},
  journal={IEEE Trans. Autom. Control.}, 
  title={Data Informativity for Robust Output Regulation}, 
  year={2024},
  volume={69},
  number={10},
  pages={7075-7080}}

@ARTICLE{9763329,
  author={Bisoffi, Andrea and De Persis, Claudio and Tesi, Pietro},
  journal={IEEE Trans. Autom. Control.}, 
  title={Controller Design for Robust Invariance From Noisy Data}, 
  year={2023},
  volume={68},
  number={1},
  pages={636-643}}

@article{Llorente_BLF,
author = {Llorente-Vidrio, Dusthon and Mera, Manuel and Salgado, Ivan and Chairez, Isaac},
title = {Robust control for state constrained systems based on composite barrier {Lyapunov} functions},
journal = {Int. J. Robust Nonlinear Control.},
volume = {30},
number = {17},
pages = {7238-7254},
year = {2020}
}

@Article{DD_Constrainted_ORP_Chinese_Paper,
author={Xia, Chaoyu
and Dong, Yi
and Wang, Chaoli
and Xu, Shengyuan},
title={Data-driven output regulation control for constrained linear systems},
journal={Science China Information Sciences},
year={2025},
month={Feb},
day={11},
volume={68},
number={3},
pages={132209},
}

@ARTICLE{10056979,
  author={dos Santos, Geovana Franca and Castelan, Eugenio B. and Lucia, Walter},
  journal={IEEE Control Syst. Lett.}, 
  title={On the Design of Constrained {PI}-Like Output Feedback Tracking Controllers via Robust Positive Invariance and Bilinear Programming}, 
  year={2023},
  volume={7},
  number={},
  pages={1429-1434}}

@INPROCEEDINGS{10644900,
  author={dos Santos, Geovana Franca and Castelan, Eugenio B. and Lucia, Walter},
  booktitle={Proc. Amer. Control Conf.}, 
  title={A Constrained Tracking Controller for Ramp and Sinusoidal Reference Signals Using Robust Positive Invariance}, 
  year={2024},
  volume={},
  number={},
  pages={1670-1675}}

@inproceedings{dos2023discrete,
  title={Discrete-time constrained {PI}-like output feedback tracking controllers-a robust positive invariance and bilinear programming approach},
  author={dos Santos, Geovana Fran{\c{c}}a and Ernesto, Jackson G and Castelan, Eug{\^e}nio B and Lucia, Walter},
  booktitle={Simp{\'o}sio Brasileiro de Automa{\c{c}}{\~a}o Inteligente-SBAI},
  volume={1},
  number={2},
  year={2023}
}

@ARTICLE{10670461,
  author={Liu, Jinhao and Yang, Jun and Yan, Yunda and Tan, Yuan and Wang, Xiangyu and Li, Shihua},
  journal={IEEE Trans. Autom. Control.}, 
  title={On Recursive Feasibility and Stability of Constrained Output Regulation}, 
  year={2025},
  volume={70},
  number={2},
  pages={1289-1296}}

@article{WANG2023111050,
title = {Constrained linear output regulation via measurement output feedback model predictive control},
journal = {Automatica},
volume = {153},
pages = {111050},
year = {2023},
issn = {0005-1098},
author = {Tianqi Wang}
}

@article{slidingmode_CRORP_SISO,
author = {Gutierrez, Ariana and Mera, Manuel and Ríos, Héctor},
title = {Integral sliding-mode-based robust output-regulation for constrained and uncertain linear systems},
journal = {Int. J. Robust Nonlinear Control.},
volume = {33},
number = {3},
pages = {2205-2218},
year = {2023}
}

@article{FRANCIS1976457,
title = {The internal model principle of control theory},
journal = {Automatica},
volume = {12},
number = {5},
pages = {457-465},
year = {1976},
issn = {0005-1098},
author = {B.A. Francis and W.M. Wonham},
}

@book{BoydLmibook,
author = {Boyd, Stephen and El Ghaoui, Laurent and Feron, Eric and Balakrishnan, Venkataramanan},
title = {Linear Matrix Inequalities in System and Control Theory},
publisher = {SIAM},
year = {1994},
address = {},
}

@ARTICLE{Blanchini_constrained_1,
  author={Blanchini, F.},
  journal={IEEE Trans. Autom. Control.}, 
  title={Feedback control for linear time-invariant systems with state and control bounds in the presence of disturbances}, 
  year={1990},
  volume={35},
  number={11},
  pages={1231-1234}}

@ARTICLE{application_spacecraft,
  author={Danielson, Claus and Kloeppel, Joseph and Petersen, Christopher},
  journal={IEEE Control Syst. Lett.}, 
  title={Spacecraft Attitude Control Using the Invariant-Set Motion-Planner}, 
  year={2022},
  volume={6},
  number={},
  pages={1700-1705}}

@ARTICLE{application_quadrotor_mobilerobot,
  author={Anevlavis, Tzanis and Liu, Zexiang and Ozay, Necmiye and Tabuada, Paulo},
  journal={IEEE Trans. Autom. Control.}, 
  title={Controlled Invariant Sets: Implicit Closed-Form Representations and Applications}, 
  year={2024},
  volume={69},
  number={7},
  pages={4506-4521}}

@book{huang2004nonlinear,
  title={Nonlinear output regulation: theory and applications},
  author={Huang, Jie},
  year={2004},
  publisher={SIAM}
}

@misc{grant2014cvx,
  title={{CVX}: {Matlab} software for disciplined convex programming, version 2.1},
  author={Grant, Michael and Boyd, Stephen},
  year={2014}
}

@article{Toh1999SDPT3,
  author    = {K. C. Toh and M. J. Todd and R. H. T{\"u}t{\"u}nc{\"u}},
  title     = {{SDPT3} --- a {MATLAB} software package for semidefinite programming},
  journal   = {Optim. Methods Softw.},
  volume    = {11},
  number    = {1--4},
  pages     = {545--581},
  year      = {1999},
  publisher = {Taylor \& Francis},
  doi       = {10.1080/10556789908805762}
}

@article{RAKOVIC2007831,
title = {Optimized robust control invariance for linear discrete-time systems: Theoretical foundations},
journal = {Automatica},
volume = {43},
number = {5},
pages = {831-841},
year = {2007},
issn = {0005-1098},
author = {S.V. Raković and E.C. Kerrigan and D.Q. Mayne and K.I. Kouramas}
}

@book{babyRudin,
  title={Principles of Mathematical Analysis},
  author={Rudin, Walter},
  year={1976},
  publisher={McGraw-Hill},
  address={New York}
}

@book{ogata2010modern,
  title={Modern Control Engineering},
  author={Ogata, Katsuhiko},
  publisher={Prentice Hall},
  address={Upper Saddle River, NJ},
  year={2010}
}

@misc{DTRCORP,
      title={Robust Cooperative Output Regulation of Discrete-Time Heterogeneous Multi-Agent Systems}, 
      author={Kursad Metehan Gul and Selahattin Burak Sarsilmaz},
      year={2026},
      eprint={2603.07783},
      archivePrefix={arXiv},
      primaryClass={eess.SY},
      url={https://arxiv.org/abs/2603.07783}, 
}

@article{koru2020cooperative,
  title={Cooperative Output Regulation of Heterogeneous Multiagent Systems: A Global Distributed Control Synthesis Approach},
  author={Koru, Ahmet Taha and Sars{\i}lmaz, Selahattin Burak and Yucelen, Tansel and Johnson, Eric},
  journal={IEEE Trans. Autom. Control},
  volume={66},
  number={9},
  pages={4289--4296},
  year={2021},
  publisher={IEEE}
}

@article{SarsilmazIJC,
author = {Selahattin Burak Sarsilmaz and Tansel Yucelen},
title = {A distributed control approach for heterogeneous linear multiagent systems},
journal = {Int. J. Control},
volume = {94},
number = {5},
pages = {1402--1414},
year = {2021},
publisher = {Taylor \& Francis}
}

@INPROCEEDINGS{8945668,
  author={Kaldmäe, Arvo and Doná, José De},
  booktitle={Proc. Aust. N. Z. Control Conf.}, 
  title={On Robustly Positively Invariant Sets and Coordinate Transformations for Discrete-time Nonlinear Systems: A Tutorial}, 
  year={2019},
  volume={},
  number={},
  pages={24-29}}

@ARTICLE{11664319,
  author={Xia, Chaoyu and Dong, Yi and Chen, Zhiyong and Chen, Jie},
  journal={IEEE Trans. Autom. Control.}, 
  title={Data-driven Output Regulation of Constrained Systems with Unmeasurable Disturbance}, 
  year={2026},
  volume={},
  number={},
  pages={1-8}}
\end{document}